\documentclass[MSNbibl,number,citesort,seceqn,dvips]{arxbj}
\usepackage{graphicx}

\aid{0}
\volume{18}
\issue{4}
\pubyear{2012}
\firstpage{1386}
\lastpage{1404}
\doi{10.3150/11-BEJ370}

\makeatletter

\newremark{example}{Example}
\newproclaim{definition}{Definition}
\newtheorem{theorem}{Theorem}[section]
\newtheorem{proposition}[theorem]{Proposition}
\newtheorem{lemma}[theorem]{Lemma}
\newremark{remark}{Remark}
\newcommand{\eqref}[1]{(\ref{#1})}
\newcommand{\ddl}{\delta}
\newcommand{\e}{\varepsilon}
\newcommand{\f}{\varphi}
\newcommand{\h}{\eta}
\newcommand{\jfd}{\psi}
\newcommand{\lft}{\lambda}
\newcommand{\m}{\mu}
\newcommand{\p}{\pi}
\newcommand{\q}{\theta}
\newcommand{\rdd}{\rho}
\newcommand{\s}{\sigma}
\newcommand{\x}{\xi}
\newcommand{\z}{\zeta}

\newcommand{\DC}{{\mathcal D}}
\newcommand{\DG}{\Delta}
\newcommand{\FC}{{\mathcal F}}
\newcommand{\LB}{{\mathbf L}}

\newcommand{\XB}{{\mathbf X}}
\newcommand{\ZB}{{\mathbf Z}}

\newcommand{\eb}{{\mathbf e}}
\newcommand{\pb}{{\mathbf p}}
\newcommand{\ub}{{\mathbf u}}
\newcommand{\wb}{{\mathbf w}}
\newcommand{\zerob}{{\mathbf0}}

\newcommand{\pbb}{{\mathbb P}}
\newcommand{\rbb}{{\mathbb R}}

\newcommand{\inv}{^{-1}}
\newcommand{\nf}{\infty}
\newcommand{\prl}{\partial}
\newcommand{\sm}{\setminus}
\newcommand{\ssdd}{\subset}

\makeatother

\begin{document}
\begin{frontmatter}

\title{Sensitivity of the limit shape of sample clouds from meta densities}

\runtitle{Sensitivity of the asymptotic behaviour of meta distributions}

\begin{aug}
\author[1]{\fnms{Guus} \snm{Balkema}\thanksref{1}\ead[label=e1]{A.A.Balkema@uva.nl}},
\author[2]{\fnms{Paul} \snm{Embrechts}\thanksref{2}\ead[label=e2]{embrechts@math.ethz.ch}}
\and
\author[3]{\fnms{Natalia} \snm{Nolde}\corref{}\thanksref{3}\ead[label=e3]{natalia@stat.ubc.ca}}

\runauthor{G. Balkema, P. Embrechts and N. Nolde}
\address[1]{Department of Mathematics, University of Amsterdam,
Science Park 904, 1098XH Amsterdam, The Netherlands. \printead{e1}}
\address[2]{Department of Mathematics and RiskLab, ETH Zurich,
Raemistrasse 101, 8092 Zurich, Switzerland. \printead{e2}}
\address[3]{Department of Statistics, University of British Columbia,
6356 Agricultural Road, Vancouver, BC V6T1Z2, Canada. \printead{e3}}
\end{aug}

\received{\smonth{12} \syear{2009}}
\revised{\smonth{12} \syear{2010}}

%
\begin{abstract}
The paper focuses on a class of light-tailed multivariate probability
distributions. These are obtained via a transformation of the margins
from a heavy-tailed original distribution. This class was introduced in
Balkema \textit{et al.}~(\textit{J. Multivariate Anal.} \textbf{101}
(2010) 1738--1754). As shown there, for the light-tailed
meta distribution the sample clouds, properly scaled, converge onto a
deterministic set. The shape of the limit set gives a good description
of the relation between extreme observations in different directions.
This paper investigates how sensitive the limit shape is to changes in
the underlying heavy-tailed distribution. Copulas fit in well with
multivariate extremes. By Galambos's theorem, existence of directional
derivatives in the upper endpoint of the copula is necessary and
sufficient for convergence of the multivariate extremes provided the
marginal maxima converge. The copula of the max-stable limit
distribution does not depend on the margins. So margins seem to play a
subsidiary role in multivariate extremes. The theory and examples
presented in this paper cast a different light on the significance of
margins. For light-tailed meta distributions, the asymptotic behaviour
is very sensitive to perturbations of the underlying heavy-tailed
original distribution, it may change drastically even when the
asymptotic behaviour of the heavy-tailed density is not affected.
\end{abstract}

%

\end{frontmatter}

\section{Introduction}


In recent years, meta distributions have been used in several
applications of multivariate probability theory, especially in finance.
The construction of meta distributions can be illustrated by a simple
example. Start with a multivariate spherical Student $t$ distribution
and transform its margins to be Gaussian. We call the new distribution
a \textit{meta distribution} with normal margins based on the
\textit{original} $t$ distribution. Since the copula of a multivariate
distribution is invariant under strictly increasing coordinatewise
transformations, the original distribution and the meta distribution
share their copula and hence have the same dependence structure. Random
vectors with Gaussian densities have asymptotically independent
components, whatever the correlation. In contrast, the above meta
distribution with Gaussian margins inherits not only the dependence
properties of the original $t$ distribution, but also the asymptotic
dependency. Asymptotic dependence is of importance in risk analysis. It
yields rank-based measures of tail dependence, and coefficients of
tail dependence (see, e.g., Chapter~5 in \cite{McNeil2005}).

Our interest is in extremes. The asymptotic behaviour of sample clouds
gives a very intuitive view of multivariate extremes. Sample clouds
from light-tailed densities tend to have a well-defined shape. The
limit shape, if it exists, describes the relation between extreme
observations in different directions; it indicates in which directions
more severe extremes are likely to occur, and how much more extreme
these will be. It has been shown in \cite{Balkema2010} that sample
clouds from meta distributions in the \textit{standard set-up}, see
below, can be scaled to converge onto a \textit{limit set}. The boundary
of this limit set has a simple explicit analytic description.
Surprisingly, the limit shape of sample clouds from the meta
distribution contains no information about the shape of the sample
clouds from the original distribution. The results of the present paper
further support this point.

Multivariate distribution functions (d.f.s) have the property that there
is a simple relation between the d.f. of the underlying random vector and
the d.f. of the coordinatewise maximum of any number of independent
observations from this distribution. One just raises the d.f. to the
given power. This makes d.f.s ideal tools to handle coordinatewise
maxima, and to study their limit behaviour. This rather analytic
approach sometimes obscures the probabilistic content of the results.
The approach via densities and probability measures on $\rbb^d$ which
is taken in this paper may at first seem clumsy, but it has the
advantage that there is a close relation to what one observes in the
sample clouds.

The aim of the paper is to investigate stability of the shape of the
limit set under changes in the original distribution. We look at
changes which do not affect the margins, or at least their asymptotic
behaviour. Keeping margins (asymptotically) unchanged allows us to
isolate the role played by the copula. We shall examine how much the
original and meta distributions in the standard set-up may be altered
without affecting the asymptotic behaviour of the scaled sample clouds.
This indicates how robust the limit shape is. Then we move on to
explore sensitivity. It turns out that the limit shape of the scaled
sample clouds from the light-tailed meta distribution is very sensitive
to certain slight perturbations of the original distribution.
Sensitivity depends on the region.

Our results suggest that the recipe outlined above for constructing
multivariate distributions with Gaussian margins with the copula of a
heavy-tailed density with a pronounced dependence structure in the
limit has to be treated with caution. The limit shape of the sample
clouds from the meta distribution is affected by perturbations of the
original heavy-tailed density, perturbations which are so small that
they do not influence the multivariate extreme value behaviour. In
going from densities with heavy-tailed margins to the meta densities
with light-tailed margins, the dependence structure of the max-stable
limit distribution is preserved by a well-known invariance result in
multivariate extreme value theory (EVT). In this paper we shall give
conditions on the severity of changes in the original heavy-tailed
distribution which are allowed if one wants to retain the asymptotic
behaviour of the coordinatewise extremes, and show that perturbations
which are negligible compared to these changes may affect the limit
shape of the sample clouds of the associated light-tailed meta
distribution.

The present paper is a follow-up to \cite{Balkema2010}. The latter
paper contains a detailed analysis of meta densities and gives the
motivation and implications of the assumptions in the standard set-up.
It presents the derivation of the analytic form in~(\ref{qE}) of the
limit shape of the sample clouds from the light-tailed meta
distribution. In the present paper, Section~\ref{sprelim} introduces
the notation and recalls the relevant definitions and results from
\cite{Balkema2010}. Section~\ref{sres} is the heart of the paper; here
we present details of the constructions which demonstrate robustness
and sensitivity of the limit shape of sample clouds from meta
distributions. Concluding remarks are given in
Section~\ref{sconc}.

\subsection*{Summary of notation}

Several of the basic symbols used in the paper are summarized in
Table~\ref{tab3}.

%
\begin{table}
\caption{Miscellaneous symbols}\label{tab3}
\begin{tabular*}{\textwidth}{@{\extracolsep{\fill}}ll@{}}
\hline
Symbol & Description \\
\hline
$f\asymp\tilde f$ & \textit{weak asymptotic equality}: ratios
$f({\mathbf x})/\tilde f({\mathbf x})$ and $\tilde f({\mathbf
x})/f({\mathbf x})$ are bounded eventually
\\
&\quad for $\|{\mathbf x}\|\to\nf$
\\
$f\sim\tilde f$ & \textit{asymptotic equality}: $\tilde f({\mathbf
x})/f({\mathbf x})\to1$ for $\|{\mathbf x}\|\to\nf$
\\
$n_D$ & the \textit{gauge function} of the set $D$: $n_D({\mathbf x})>0$
for ${\mathbf x}\neq\zerob$, $D=\{n_D<1\}$ and
\\
&\quad $n_D(c{\mathbf x})=cn_D({\mathbf x})$ for ${\mathbf x}\in\rbb^d$, $c\ge0$\\
$\eb$; $E_\times$ & a vector of ones in $\rbb^d$; the diagonal cross
(see~(\ref{qrsdx}))
\\
$B$, $C$ & the open Euclidean unit ball in $\rbb^d$; the open cube
$(-1,1)^d$
\\\hline
\end{tabular*}
\end{table}

Throughout the paper, it is convenient to keep in mind two spaces:
${\mathbf z}$-space on which the heavy-tailed d.f.s $F,F^*, \ldots$
are defined,
and ${\mathbf x}$-space on which the light-tailed meta d.f.s
$G,G^*,\ldots$ are
defined. Table~\ref{tab1} compares notation used for mathematical
objects on these two spaces.

%
\begin{table}
\tabcolsep=3pt
\caption{Symbols used to distinguish various objects of interest in
${\mathbf z}$-space and in ${\mathbf x}$-space. Notation for margins
assumes that all
marginal densities are equal and symmetric}\label{tab1}
\begin{tabular*}{\textwidth}{@{\extracolsep{\fill}}lll@{}}
\hline
${\mathbf z}$-space & ${\mathbf x}$-space & Comments \\
\hline
$F$ and $f$ & $G=F\circ K$ and $g$ & joint d.f. and density \\
$F_0$ and $f_0$ & $G_0=F_0\circ K_0$ and $g_0$ & margin d.f. and density
\\
$\p$ & $\m$ & probability measures \\
$\p_n$ & $\m_n$ & mean measures of scaled sample clouds\\
$\ZB$, $\ZB_1,\ZB_2,\ldots$ & $\XB$, $\XB_1,\XB_2,\ldots$ &
random vectors \\
$N_n$ & $M_n$ & scaled $n$-point sample clouds \\
$N$: Poisson point process & $E$: limit set & limit of scaled sample
clouds \\
$c_n\dvt1-F_0(c_n)\sim1/n$ & $b_n\dvt-\log(1-G_0(b_n))\sim\log n$ &
scaling constants \\
$(B_n)$, $t_n=K_0(s_n)$ & $(A_n)$, $s_n$ & block partitions, division
points \\ \hline
\end{tabular*}
\end{table}

\section{Preliminaries}\label{sprelim}

\subsection{Definitions and standard set-up}

Meta distributions are constructed by transforming margins of a given
multivariate distribution so as to obtain a new (meta) distribution
with desired margins. To fix our notation and terminology, let us
present this construction procedure formally.

Consider a random vector $\ZB$ in $\rbb^d$ with d.f. $F$ and continuous
margins $F_i$, $i=1,\ldots,d$. Let $G_1,\ldots,G_d$ be continuous d.f.s
on $\rbb$ which are strictly increasing on the intervals
$I_i=\{0<G_i<1\}$. Define the transformation
\begin{equation}\label{eq:Ki}
K(x_1,\ldots,x_d)=(K_1(x_1),\ldots,K_d(x_d)),\qquad
K_i(s)=F_i^{-1}(G_i(s)),\qquad i=1,\ldots,d.
\end{equation}
The d.f. $G=F\circ K$ is
called the \textit{meta distribution} with \textit{margins} $G_i$ based on
the \textit{original} d.f. $F$. The coordinatewise map
$K=K_1\otimes\cdots\otimes K_d$, which maps ${\mathbf x}=(x_1,\ldots
,x_d)\in
I=I_1\times\cdots\times I_d$ into the vector
${\mathbf z}=(K_1(x_1),\ldots,K_d(x_d))$, is called the \textit{meta
transformation}.


In the paper, we restrict attention to meta distributions with
light-tailed margins based on a heavy-tailed original distribution. The
precise assumptions we make on the meta margins and the original
distribution are summarized in the standard set-up given in
Definition~\ref{dssu} below. But before stating assumptions of the
standard set-up, we need several additional definitions and notation.

Recall the basic example we started with in the introduction of a meta
distribution with normal margins based on a multivariate $t$
distribution. The multivariate $t$ density has a simple structure. It
is fully characterized by the shape of its level sets, scaled copies of
the defining ellipsoid, and by the decay $c/r^{\lft+d}$ of its tails
along rays. The constant $\lft>0$ denotes the degrees of freedom, $d$ the
dimension of the underlying space, and $c$ is a positive constant
depending on the direction of the ray. In the more general setting of
the paper, the tails of the density are allowed to decrease as
$cL(r)/r^{\lft+d}$ for some slowly varying function $L$ and the condition
of elliptical level sets is replaced by the requirement that the level
sets are equal to scaled copies of a fixed bounded convex or
star-shaped set (a set $D$ is star-shaped if ${\mathbf z}\in D$ implies
$t{\mathbf z}\in D$ for $0\le t<1$). Due to the power decay of the
tails, the
density $f$ is said to be \textit{heavy-tailed}. Densities with the above
properties constitute the class $\FC_\lft$.

\begin{definition}\label{dFC}
The set $\FC_\lft$ for $\lft>0$ consists of all positive
continuous densities $f$ on $\rbb^d$ which are asymptotic to a function
of the form $f_*(n_D({\mathbf z}))$, where $f_*(r)=L(r)/r^{(\lft+d)}$
is a
continuous decreasing function on $[0,\nf)$, $L$ varies slowly, and
$n_D$ is the gauge function of the set $D$ (see Table~\textup{\ref
{tab3}}). The
set $D$ is bounded, open and star-shaped. It contains the origin and
the gauge function is continuous.
\end{definition}

The reader may keep in mind the case where $D$ is a convex symmetric
set, $D=-D$. In that case, the gauge function is a norm, and $D$ the
open unit ball in this norm. In the standard set-up, the normal margins
of the meta density in the introduction are generalized to include
densities whose tails are asymptotic to a \textit{von Mises function}:
$g_0(s)\sim\mathrm{e}^{-\jfd(s)}$ for $s\to\nf$ with \textit{scale function}
$a=1/\jfd'$, where $\jfd$ is a $C^2$ function with a positive derivative
such that
\begin{equation}\label{qvMfn}
\jfd(s)\to\nf,\qquad a(s)'\to0,\qquad s\to\nf.
\end{equation}
This condition on the meta margins ensures that they lie in the maximum
domain of attraction of the Gumbel limit law $\exp(-\mathrm{e}^{-x})$,
$x\in\rbb$; see, for example, Proposition~1.4 in \cite{Resnick1987}.

We are now able to define the standard set-up.
\begin{definition}\label{dssu}
Let $f$ be a density in $\FC_\lft$ for some $\lft
>0$ with
margins all equal to a positive continuous symmetric density $f_0$, and
let $g_0$ be a continuous, positive, symmetric density on $\rbb$
asymptotic to a von Mises function $\mathrm{e}^{-\jfd}$, where in addition
to~\textup{(\ref{qvMfn})} $\jfd$ is assumed to vary regularly at
infinity with
exponent $\q>0$. In the \textit{standard set-up}, the meta density $g$ is
based on $f$ and has margins $g_0$.
\end{definition}

This set-up is the same as in~\cite{Balkema2010} except that for
simplicity we additionally assume that the original density has equal
margins. As a consequence, the components of the meta transformation
$K$ are equal:
\begin{equation}\label{qemK}
K\dvtx{\mathbf x}\mapsto{\mathbf z}=(K_0(x_1),\ldots
,K_0(x_d)),\qquad
K_0=F_0\inv\circ G_0,\qquad K_0(-t)=-K_0(t).
\end{equation}

\subsection{Convergence of sample clouds}

An \textit{$n$-point sample cloud} is the point process consisting of the
first $n$ points of a sequence of independent observations from a given
distribution, after proper scaling. We write
\begin{equation}\label{q1Nn}
N_n=\{\ZB_1/a_n,\ldots,\ZB_n/a_n\},
\end{equation}
where $\ZB_1,\ZB_2,\ldots
$ are
independent observations from the given probability distribution on
$\rbb^d$, and $a_n$ are positive scaling constants. It is customary to
write $N_n(A)$ for the number of the points of the sample cloud that
fall into the set $A$.

In this section, we discuss the asymptotic behaviour of sample clouds
from the original density $f$ and from the associated meta density $g$
in the standard set-up. The difference in the asymptotic behaviour is
striking: sample clouds from the heavy-tailed density $f$ converge in
distribution to a Poisson point process on $\rbb^d\sm\{\zerob\}$
whereas sample clouds from the light-tailed meta density $g$ tend to
have a clearly defined boundary. They converge onto a deterministic
set.

\subsubsection{\texorpdfstring{Convergence for densities in $\FC_\lft$ and measures in $\DC_\lft$}
{Convergence for densities in F lambda and measures in D lambda}}

For densities in $\FC_\lft$, $\lft>0$, there is a simple limit relation
which follows from the regular variation of the function $f_*$ in
Definition~\ref{dFC} and from\vadjust{\goodbreak} the homogeneity of the gauge function:
\begin{equation}\label{qsls0}
\dfrac{f(r_n\wb_n)}{f_*(r_n)}\to h(\wb),\qquad
\wb_n\to\wb\neq\zerob,r_n\to\nf,
\end{equation}
where
\begin{equation}\label{q1h}
h(\wb)=1/n_D(\wb)^{\lft+d}.
\end{equation}
Convergence is uniform and $\LB^1$ on
the complement of centered balls. If $\ZB_1,\ZB_2,\ldots$ are
independent observations from the density $f$ then the sample clouds
$N_n$ in~(\ref{q1Nn}) converge in distribution to the Poisson point
process with intensity $h$ weakly on the complement of centered balls
for a suitable choice of scaling constants $a_n$.

A probability measure $\p$ on $\rbb^d$ varies regularly with scaling
function $a(t)\to\nf$ if there is an infinite Radon measure $\rdd$ on
$\rbb^d\sm\{\zerob\}$ such that the finite measures $t\p$ scaled by
$a(t)$ converge to $\rdd$ vaguely on $\rbb^d\sm\{\zerob\}$. The measure
$\rdd$ has the scaling property
\begin{equation}\label{q1r}
\rdd(rA)=\rdd(A)/r^\lft,\qquad r>0
\end{equation}
for all Borel sets $A$ in $\rbb^d\sm\{\zerob\}$. The constant
$\lft\ge0$ is the exponent of regular variation. If it is positive, then
$\rdd$ gives finite mass to the complement of the open unit ball $B$ and
weak convergence holds on the complement of centered balls. We shall
denote the set of all probability measures which vary regularly with
exponent $\lft>0$ and with limit measure $\rdd$ by $\DC_\lft(\rdd
)$ (or just
$\DC_\lft$). In particular, $\p\in\DC_\lft$ if $\p$ has a continuous
density in $\FC_{\lft}$. As above, for independent observations
$\ZB_1,\ZB_2,\ldots$ from a distribution $\p\in\DC_\lft(\rdd)$
the scaled
sample clouds $N_n$ in~(\ref{q1Nn}) (with scaling constants $a_n$)
converge in distribution to a Poisson point process with mean measure
$\rdd$ weakly on the complement of centered balls (since the mean
measures converge; see, e.g., Proposition~3.21 in \cite{Resnick1987} or
Theorem~11.2.V in \cite{Daley2008}).

\subsubsection{Convergence for meta densities in the standard set-up}

Sample clouds from light-tailed meta densities in the standard set-up,
under suitable scaling, converge onto a deterministic set, referred to
as the \textit{limit set}, in the sense of the following definition.

\begin{definition}\label{dconv}
For a compact set $E$ in $\rbb^d$, the $n$-point
point processes $N_n$ \textit{converge onto} $E$ as $n\to\nf$ if for open
sets $U$ containing $E$, the probability of a point outside $U$
vanishes, $\pbb\{N_n(U^c)>0\}\to0$, and if
\[
\pbb\{N_n(\pb+\e B)>m\}\to1\qquad\mbox{for all $m\ge1$ and $\e
>0$, $\pb\in E$.}
\]
\end{definition}

We now recall Theorem~2.6 of \cite{Balkema2010}, which characterizes
the shape of the limit set for meta distributions in the standard
set-up.

\begin{theorem}\label{tE}
Let the meta density $g$ satisfy the assumptions of the
standard set-up of Definition~\textup{\ref{dssu}}. Define
\begin{equation}\label{qE}
E:=E_{\lft,\q}=\{\ub\in\rbb^d\mid|u_1|^\q+\cdots+|u_d|^\q+\lft
\ge
(\lft+d)\|\ub\|_\nf^\q\}.
\end{equation}
If $\jfd(r_n)\sim\log n$, then for the sequence of independent
observations $\XB_1,\XB_2,\ldots$ from the meta density $g$ the sample
clouds $M_n=\{\XB_1/r_n,\ldots,\XB_n/r_n\}$ converge onto $E$.
\end{theorem}

\begin{figure}
\begin{tabular}{@{}cc@{}}

\includegraphics{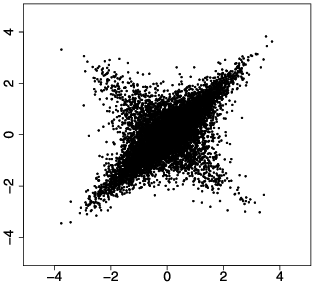}
&\includegraphics{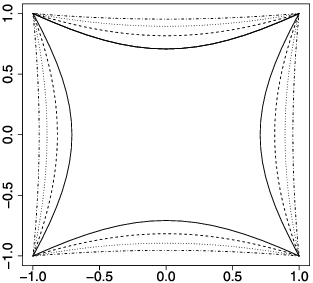}\\
(a)&(b)\\
\end{tabular}
\caption{Panel~(a): A sample cloud of 10,000 points from a meta
distribution with standard Gaussian margins based on a Cauchy
distribution. Panel~(b): The limit sets $E_{\lft,\q}$ in \protect\eqref{qE} with
$\q=2$ and $\lft=1,2,4,10$ (solid, dash, dot, dotdash lines,
resp.).}
\label{f1} 
\end{figure}

Figure~\ref{f1} illustrates the limit set $E_{\lft,\q}$ in the above
theorem for several values of the parameters $\lft$ and $\q$. The sample
cloud in Panel~(a) of the figure is simulated from the meta
distribution with standard Gaussian margins based on a Cauchy
distribution with elliptic level sets. Hence, the associated limit set
has parameters $\lft=1$ and $\q=2$.

%

\subsection{Further notation and conventions}

In order to ease the exposition, we introduce some additional
assumptions and notation.

All univariate d.f.s are assumed to be continuous and strictly
increasing. The d.f.s $F_0$ and $\tilde F_0$ on $\rbb$ are \textit{tail
asymptotic} if
\[
\tilde F_0(-t)/F_0(-t)\to1,\qquad\bigl(1-\tilde F_0(t)\bigr)/\bigl(1-F_0(t)\bigr)\to
1,\qquad t\to\nf.
\]
The sample clouds from a heavy-tailed d.f. $\tilde F$ \textit{converge to
the point process} $\tilde N$ if $\tilde N$ is a Poisson point process
on $\rbb^d\sm\{\zerob\}$, and if the sample clouds converge to
$\tilde
N$ in distribution weakly on the complement of centered balls, where
the scaling constants $c_n$ satisfy $1-\tilde F_0(c_n)\sim1/n$. Two
heavy-tailed d.f.s $F^*$ and $F^{**}$ \textit{have the same asymptotics} if
the margins are tail asymptotic and if the sample clouds converge to
the same point process. The light-tailed d.f.s $G^*$ and $G^{**}$
\textit{have the same asymptotics} if the margins are tail asymptotic and
the sample clouds converge onto the same compact set $E^*$, with scale
factors $b_n$ which satisfy $-\log(1-G_0(b_n))\sim\log n$. One may
replace a scaling sequence by a sequence asymptotic to it without
affecting the limit. The scaling of sample clouds is determined up to
asymptotic equality by the margins since the univariate projections
have to converge. Tail asymptotic margins yield asymptotic scalings.

\section{Results}\label{sres}

\subsection{Preview}

We now turn our attention to the main issue of the paper. We start out
with a pair of d.f.s $F$ and $G=F\circ K$ which satisfy the conditions of
the standard set-up. Introduce a new pair $(F^*,G^*)$ where $F^*$ and
$G^*=F^*\circ K$ are related by the same meta transformation $K$, and
where one term of the pair has the same asymptotic behaviour as the
corresponding term in $(F,G)$. What does this imply for the other term
of the pair?
\begin{enumerate}[(Q.2)]
\item[(Q.1)] If the scaled sample clouds from $G^*$ and from $G$
converge onto the same set $E$, do the scaled sample clouds from $F^*$
converge to the same point process $N$ as those from $F$?
\item[(Q.2)]
If the scaled sample clouds from $F^*$ and from $F$ converge to the
same Poisson point process $N$, do the scaled sample clouds from $G^*$
converge onto the same set $E$ as those from $G$?
\end{enumerate}

For coordinatewise maxima, convergence is determined by the copula and
the answer to the questions is ``Yes.'' If the margins lie in the
univariate domain of attraction and $F$ and $G$ have the same copula
then convergence of the coordinatewise maxima and the associated sample
clouds for one d.f. implies convergence for the other d.f. In the
situation of this paper the answer to both questions is ``No.'' For
example, for (Q.1), if we replace $f$ by a weakly asymptotic density
$f^*\asymp f$, the asymptotic behaviour of sample clouds from $g^*$ is
not affected, since $g^*\asymp g$ (see Lemma~\ref{mrs0} below), but the
scaled sample clouds from $f^*$ obviously need not converge. This
section contains some counterexamples which will be worked out further
in the next two sections.

What we want to do is specify the margins $f_0$ and $g_0$ (which
determine the meta transformation $K$), then vary the copula and check
the limit behaviour of the sample clouds (where we impose the condition
that both converge). We are looking for d.f.s $F^*$ and $G^*$ with the
properties:
\begin{enumerate}[(P.3)]
\item[(P.1)] $F^*$ has marginal densities $f_0$;
\item[(P.2)] $G^*$ is the meta distribution based on $F^*$ with
marginal densities $g_0$;
\item[(P.3)] the scaled sample clouds from $F^*$ converge to a Poisson
point process $N^*$;
\item[(P.4)] the scaled sample clouds from $G^*$ converge onto a
compact set $E^*$.
\end{enumerate}

Moreover, we would like $E^*$ to be the set $E_{\lft,\q}$ in \eqref{qE},
or $N^*$ to have mean measure $\rdd^*=\rdd$ with intensity $h$ in
\eqref{q1h}. So we either choose $F^*$ to have the same asymptotics as
$F$, or $G^*$ to have the same asymptotics as $G$. Note that the
conditions (P1)--(P4) have certain implications. The mean measure
$\rdd^*$ of the Poisson point process $N^*$ is an excess measure with
exponent $\lft$, see~(\ref{q1r}), its margins are equal and symmetric
with intensity $\lft/|t|^{\lft+1}$ since the marginal densities $f_0$ are
equal and symmetric and the scaling constants $c_n$ ensure that
$\rdd\{w_d\ge1\}=1$. The limit set $E^*$ is a subset of the cube
$C=[-1,1]^d$ and projects onto the interval $[-1,1]$ in each
coordinate, again by our choice of the scaling constants.

The two sections below describe procedures for altering distributions
without changing the margins too much. The first procedure uses block
partitions. A\vadjust{\goodbreak} block partition is a special kind of partition into
coordinate blocks. If the blocks are relatively small, then the
asymptotics of a distribution do not change if one replaces it by one
which gives the same mass to each block. Block partitions are mapped
into block partitions by the meta transformation $K$. The mass is
preserved, but the size and shape of the blocks may change drastically.
The block partitions provide insight in the relation between the
asymptotic behaviour of the measures $\mathrm{d}F^*$ and $\mathrm{d}G^*$. In the second
procedure, we replace $\mathrm{d}F$ by a probability measure $\mathrm{d}\tilde F$, which
agrees outside a bounded set with a mixture $\mathrm{d}(F+F^o)$, where $F^o$ has
lighter margins than $F$:
\begin{equation}\label{qrsneg}
F_j^o(-t)\ll F_0(-t),\qquad 1-F_j^o(t)\ll 1-F_0(t),\qquad t\to\nf,\
j=1,\ldots,d.
\end{equation}
\noindent This condition ensures that $\tilde F$ and $F$ have
the same asymptotics. The two corresponding light-tailed meta d.f.s
$\tilde G$ and $G$ on ${\mathbf x}$-space may have different
asymptotics since
the scaling constants $b^o_n$ and $b_n$ may be asymptotic even though
$G^o$ has lighter tails than $G$. If this is the case, and the scaled
sample clouds from $G^o$ converge onto a compact set $E^o$, then those
from $\tilde G$ converge onto the union $E\cup E^o$, which may be
larger than $E$. These two procedures enable us to construct d.f.s $F^*$
and $G^*=F^*\circ K$ with the same marginal densities as the d.f.s $F$
and $G$ in the standard set-up which exhibit unexpected behaviour:
\begin{enumerate}[(Ex.2)]
\item[(Ex.1)] $G^*$ and $G$ have the same asymptotics, but the scaled
sample clouds from $F^*$ converge to a Poisson point process which
lives on the diagonal. (See Proposition~\ref{thc1} and
Example~\ref{em1}.)
\item[(Ex.2)] The scaled sample clouds from $G^*$ converge onto
$E^*=A\cup E_{\times}$, where $E_{\times}$ is the \textit{diagonal cross}
\begin{equation}\label{qrsdx} E_{\times}=\{r\ddl\mid0\le r\le1,
\ddl\in\{-1,1\}^d\},
\end{equation}
and
$A$ is the closure of an open star-shaped set $D\ssdd(-1,1)^d$ with a
continuous gauge function $n_D$, see Table~\ref{tab3} above. The d.f.s $F^*$ and
$F$ have the same asymptotics. The density $f^*$ is asymptotic to $f$
on every ray. (See Proposition~\ref{thmix} and Example~\ref{em1}.)
\end{enumerate}

What does the copula say about the asymptotics? Everything, since it
determines the d.f. if the margins are given; nothing, since the examples
above show that there is no relation between the asymptotics of $F^*$
and the asymptotics of $G^*$ even with the prescribed margins $f_0$ and
$g_0$. One might hope that at least the parameters $\lft$ and $\q$,
determined by the margins, might be preserved in the asymptotics. The
point process $N^*$ reflects the parameter $\lft$ in the marginal
intensities $\lft/|t|^{\lft+1}$. However, $E^*=E_{\lft^*,\q^*}$ may
hold for
any $\lft^*$ and $\q^*$ in $(0,\nf)$ by taking $A=E_{\lft^*,\q^*}$ in
(Ex.2) above.

We now start with the technical details.

The construction procedures discussed above will change an original d.f.
$\hat F$ with margins $F_0$ into a d.f. $\tilde F$ whose margins
$\tilde
F_j$ are tail equivalent to $F_0$. This is no serious obstacle.

\begin{proposition}\label{prs0}
Let the scaled sample clouds from
$\tilde F$ converge
to a point process $\tilde N$, and let the scaled sample clouds from
$\tilde G=\tilde F\circ K$ converge onto the compact set $\tilde E$. If
the margins $\tilde F_j$ are continuous and strictly increasing and
tail asymptotic to $F_0$, then there exists a d.f. $F^*$ with margins
$F_0$ such that $F^*$ has the same asymptotics as $\tilde F$ and
$G^*=F^*\circ K$ has the same asymptotics as $\tilde G$.

If moreover $\tilde F$ has a density $\tilde f$ with margins asymptotic
to $f_0$ at $\pm\nf$, then $F^*$ has a density $f^*$, and for any
vector $\wb$ with nonzero coordinates and any sequence $\wb_n\to\wb$
and $r_n\to\nf$ there is a sequence $\wb'_n\to\wb$ such that
$f^*(r_n\wb_n)\sim\tilde f(r_n\wb'_n)$.
\end{proposition}

\begin{pf}
Let $F^*=\tilde
F\circ K_F$ be the meta d.f. based on $\tilde F$ with margins $F_0$. The
components $K_{Fj}=\tilde F_j\inv\circ F_0$ are homeomorphisms and
satisfy $K_{Fj}(t)\sim t$ for $|t|\to\nf$. (Here, we use that the
marginal tails vary regularly with exponent $-\lft\ne0$.) It follows that
\begin{equation}\label{qrskf}
 \|K_F({\mathbf z})-{\mathbf z}\|/\|
{\mathbf z}\|\to
0,\qquad\|{\mathbf z}\|\to\nf.
\end{equation}
 This
ensures that $\tilde F$ and $F^*$ have the same asymptotics. (For any
$\e>0$, the distance between the scaled sample point $\ZB/c_n$ and
$K_F(\ZB)/c_n$ is bounded by $\e\|\ZB\|/c_n$ for $\|\ZB\|\ge\e
c_n$ and
$n\ge n_\e$.) A similar argument shows that $\tilde G=\tilde F\circ K$
and $G^*=F^*\circ K$, the meta d.f. based on $\tilde G$ with margins
$G_0$, have the same asymptotics. Here, we use that $\tilde G_j=\tilde
F_j\circ K_0$ is tail asymptotic to $G_0=F_0\circ K_0$ since $\tilde
F_j$ is tail asymptotic to $F_0$. Under the assumptions on the density
the Jacobian of $K_F$ is asymptotic to one in the points $r_n\wb'_n$
and~(\ref{qrskf}) gives the limit relation with
$r_n\wb'_n=K_F\inv(r_n\wb_n)$.
\end{pf}

In general, the densities $f^*$ and $\tilde f$ (in the notation of the
above proposition) are only weakly asymptotic, as in Proposition~1.8 in
\cite{Balkema2010}. The density $f^*$ on ${\mathbf z}$-space is
related to $f$
in the same way as the density $g^*$ is related to $g$. If $f^*\asymp
f$ or $f^*\le Cf$ or $f^*\sim f$, then these relations also holds for
$g^*$ and $g$, and vice versa. Similarly, for the margins: $g^*_j\sim
g_0$ at $\nf$ implies $f^*_j\sim f_0$ at $\nf$. These results are
formalized in the lemma below.

\begin{lemma}\label{mrs0} If $F^*$ has density $f^*$ with margins
$f_j^*$ and
$G^*=F^*\circ K$ has density $g^*$ with margins $g^*_j$, then
\[
g^*({\mathbf x})/g({\mathbf x})=f^*({\mathbf z})/f({\mathbf z}),\qquad
g^*_j(s)/g_0(s)=f^*_j(t)/f_0(t),\qquad{\mathbf z}=K({\mathbf
x}),\qquad t=K_0(s).
\]
\end{lemma}

\begin{pf}
The Jacobian drops out in the quotients.
\end{pf}

\subsection{Block partitions}\label{sdom2}

Block partitions are a good tool for testing the robustness and
sensitivity of asymptotic behaviour of sample clouds via simple
constructions.

Consider a partition of $\rbb^d$ into bounded Borel sets $A_n$ given by
coordinate blocks. Since our d.f.s have continuous margins, the
boundaries of the blocks are null sets, and we shall not bother about
boundary points, and treat the blocks as open sets. To construct such a
\textit{block partition} start with an increasing sequence of cubes
\[
s_nC=(-s_n,s_n)^d,\qquad 0<s_1<s_2<\cdots, s_n\to\nf.
\]
Subdivide the ring $R_n=s_{n+1}C\sm s_nC$ between two successive cubes
into blocks by a symmetric partition of the interval
$[-s_{n+1},s_{n+1}]$ with division points $\pm s_{nj}$,
$j=1,\ldots,m_n$ with $s_{nm_n}=s_n$, $s_{nm_n+1}=s_{n+1}$ and
\[
-s_{n+1}<-s_n<\cdots<-s_{n1}<s_{n1}<\cdots<s_n<s_{n+1}.
\]
This gives a partition of the cube $s_{n+1}C$ into $(2m_n+1)^d$ blocks
of which $(2m_n-1)^d$ form the cube $s_nC$. The remaining blocks form
the ring $R_n$. The meta transformation $K$ transforms block partitions
in ${\mathbf x}$-space into block partitions in ${\mathbf z}$-space. A
comparison of
the original block partition with its transform gives a good indication
of the way in which the meta transformation distorts space.

The block partition introduced above is \textit{regular} if $s_{n+1}\sim
s_n$ and $\Delta_n/s_n\to0$, where $\Delta_n$ is the maximum of
$s_{n1}, s_{n2}-s_{n1},\ldots,s_{nm_n}-s_{nm_n-1}$.

Our first aim is to investigate how much the d.f.s $F$ and $G=F\circ K$
in the standard set-up may be altered without affecting the asymptotic
behaviour of the scaled sample clouds. Regular partitions give a simple
answer to a related question: If one replaces density $f$ or $g$ by a
discrete distribution, how far apart are the atoms allowed to be if one
wants to retain the asymptotic behaviour of the sample clouds from the
given density?

\begin{proposition}\label{prs1}
Let $A_1,A_2,\ldots$ in ${\mathbf x}$-space be a
regular block
partition. Suppose the sample clouds from the probability distribution
$\m$ scaled by $r_n$ converge onto the compact set $E$. Let $\tilde\m$
be a probability measure such that $\tilde\m(A_n)=\m(A_n)$ for $n\ge
n_0$. Then the sample clouds from $\tilde\m$ scaled by $r_n$ converge
onto $E$.
\end{proposition}

\begin{pf}
Let $\pb\in E$, and $\e>0$. Let $\m_n$ denote the
mean measure from the scaled sample cloud from $\m$ and $\tilde\m_n$
the same for $\tilde\m$. Then $\m_n(\pb+(\e/2)B)\to\nf$. Because the
sets $A_n$ are relatively small there exists $n_1$ such that any set
$A_n$ which intersects the ball $r_n\pb+r_n\e B$ with $n\ge n_1$ has
diameter less than $\e r_n/2$. Let $U_n$ be the union of the sets $A_n$
which intersect $r_n\pb+(r_n\e/2)B$. Then $U_n\ssdd r_n\pb+\e r_nB$ and
hence
\[
\m_n\bigl(\pb+(\e/2)B\bigr)\le\m_n(U_n/r_n)=\tilde\m_n(U_n/r_n)\le\tilde
\m_n(\pb+\e B).
\]
Similarly, $\tilde\m_n(U^c)\to0$ for any open set $U$ which contains
$E$.
\end{pf}

\begin{remark}
The result also holds if $\tilde\m(A_n)\sim\m
(A_n)$ provided
$\m(A_n)$ is positive eventually.
\end{remark}

There is an analogous result for regular partitions in ${\mathbf z}$-space.

\begin{proposition}\label{prs2}
Suppose $\p\in\DC_\lft(\rdd)$
with scaling constants
$c_n$. Let $A_1,A_2,\ldots$ be a regular block partition and let
$\tilde\p$ be a probability measure on $\rbb^d$ such that
$\tilde\p(A_n)=\p(A_n)$ for $n\ge n_0$. Then $\tilde\p\in\DC
_\lft
(\rdd)$
with scaling constants $c_n$.
\end{proposition}

\begin{pf}Any closed block
$A\ssdd\rbb^d\sm\{\zerob\}$ whose boundary carries no $\rdd$-mass is
contained in an open block $U$ with $\rdd(U)<\rdd(A)+\e$. As in the proof
of the previous proposition, for $n\ge n_1$ there is a union $U_n$ of
atoms $A_n$ such that $A\ssdd U_n/c_n\ssdd U$. By homogeneity of $\rdd
$, the
boundary of a block has positive mass $\rdd(\prl B_n)>0$ only if one of
the vertices lies in a coordinate plane. By assumption, the division
points are nonzero.
\end{pf}

We thus have the following simple situation: $A_1,A_2,\ldots$ is a
block partition in ${\mathbf x}$-space and
$B_1=K(A_1),B_2=K(A_2),\ldots$ the
corresponding block partition in ${\mathbf z}$-space, where $K$ is the meta
transformation in \eqref{qemK}. Let $\tilde\p$ be a probability measure
in ${\mathbf z}$-space and $\tilde\m$ a probability measure in
${\mathbf x}$-space,
linked by $K$, that is, $\tilde\p=K(\tilde\m)$. Then
$\tilde\p(A_n)=\tilde\m(B_n)$ for all $n$. So
\begin{equation}\label{qrs1}
\tilde\p(A_n)\sim\int_{A_n}f({\mathbf z})\,\mathrm{d}{\mathbf z}\quad\iff\quad
\tilde\m(B_n)\sim\int_{B_n}g({\mathbf x})\,\mathrm{d}{\mathbf x}.
\end{equation}

\begin{theorem}\label{thrs3} Assume the standard set-up in
Definition~\textup{\ref{dssu}}.
If the two block partitions $A_n$ and $B_n=K(A_n)$ both are regular,
and one of the equivalent asymptotic equalities in~\eqref{qrs1} holds,
then the sample clouds from $\tilde\p$ scaled by $c_n$ converge to the
Poisson point process with intensity $h$ in \eqref{q1h}, and the sample
clouds from $\tilde\m$ scaled by $r_n$ converge onto the set
$E=E_{\lft,\q}$ in \eqref{qE}.
\end{theorem}

\begin{pf}Combine Propositions~\ref{prs1}
and~\ref{prs2}.
\end{pf}

Unfortunately, the meta transformation $K$ is very nonlinear.
Regularity of one block partition does not imply regularity of the
other block partition.

We first consider the case when the block partition $(A_n)$ in
${\mathbf x}$-space is regular, but $(B_n)$ is not. The block partition $(A_n)$
in ${\mathbf x}$-space is based on a sequence of cubes $s_nC=(-s_n,s_n)^d$.
Successive cubes are of the same size asymptotically, $s_{n+1}\sim
s_n$. The cubes $t_nC$ in ${\mathbf z}$-space with $t_n=K_0(s_n)$ may
grow very
fast. It is possible that $t_n\ll t_{n+1}$. The corresponding partition
with blocks $B_n=K(A_n)$ in ${\mathbf z}$-space then certainly is not regular.

\begin{proposition}\label{prs4} Assume the standard set-up. Let $\h
\in(0,1)$.
There is
a sequence $0<s_1<s_2<\cdots$ such that $s_n\to\nf$ and $s_{n+1}\sim
s_n$, and such that
\[
t_n=K_0(s_n)=n^{n^{n^\h}}.
\]
\end{proposition}
\begin{pf}
We have $g_0\sim \mathrm{e}^{-\jfd}$. This implies
$1-G_0(s)\sim
a(s)g_0(s)\sim \mathrm{e}^{-\Psi(s)}$, where $\Psi$ like $\jfd$ varies regularly
with exponent $\q$. Write $s_n=\mathrm{e}^{\s_n}$ and $\tau\Psi(s_n)\sim
\mathrm{e}^{\q
r(\s_n)}$, where $r$ is a $C^2$ function with $r'(t)\to1$ and
$r''(t)\to0$ as $t\to\nf$, and $\tau:=1/\lft$. It has been shown in
\cite{Balkema2010} (Equation~(1.13)) that
\[
K_0(s)=t\sim c\mathrm{e}^{\f(s)},\qquad\f(s)=\tau q(\Psi(s))\sim\tau\Psi
(s),\qquad s\to\nf
\]
for some positive constant $c$. This gives $\log t_n=\log K_0(s_n)\sim
\mathrm{e}^{\q r(\s_n)}$. Since $\log\log t_n=n^{1-\e}\log n+\log\log n$ has
increments which go to zero, so does $\q r(\s_n)$, and hence also
$\s_n$ since $r'$ tends to one. It follows that $s_{n+1}\sim s_n$.
\end{pf}

Choose $s_{nm_n-1}=s_{n-1}$. Then the cube $(s_{n-1}\eb,s_{n+1}\eb)$ is
a union of $2^d$ blocks in the partition on ${\mathbf x}$-space, and so
is the
cube $(t_{n-1}\eb,t_{n+1}\eb)$ in ${\mathbf z}$-space; $\eb
=(1,\ldots,1)$
denotes a vector of ones in $\rbb^d$. The union $U$ of these latter
cubes has the property that the scaled sets $U/t_n$ converge to
$(0,\nf)^d$ for $t_n\to\nf$ if $t_n\ll t_{n+1}$.

%
%
\begin{figure}
\begin{tabular}{@{}cc@{}}

\includegraphics{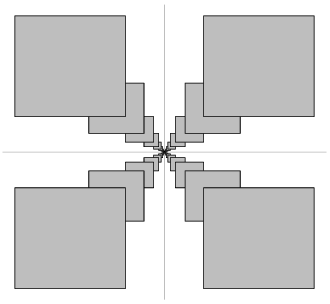}
&\includegraphics{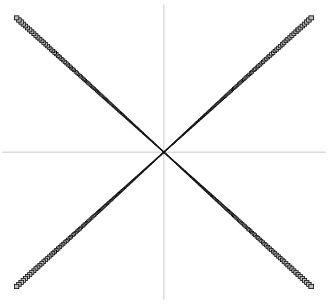}\\
(a)&(b)\\
\end{tabular}
\caption{Panel~(a): Schematic representation of the four sets
$U_{\Delta}=\bigcup_n [\Delta_1 t_{n-1}\eb,\Delta_2 t_{n+1}\eb]$ in
$\rbb^2$ for $\Delta=(\Delta_1,\Delta_2)\in\{-1,1\}^2$; see
Proposition~\protect\ref{thc1}. Panel~(b): the image of the sets $U_{\Delta}$
in ${\mathbf x}$-space. The division points of the partitions $(A_n)$ and
$(B_n)$ are $s_n=n^{\sqrt{n}}\log n$ and
$t_n=\mathrm{e}^{s_n}=n^{n^{\sqrt{n}}}$, respectively. Use
Proposition~\protect\ref{prs4} with $K_0(s)=\mathrm{e}^s$, $s>s_0$ and $\e=0.5$. (a) \textbf{z}-space. (b) \textbf{x}-space (regular).}
\label{fcase1} 
\end{figure}
%

\begin{proposition}\label{thc1} Assume the standard set-up. The excess
measure $\rdd$
has the continuous density $h$ in \textup{(2.6)} and does not charge the
coordinate planes. Moreover, $\rdd\{\ddl z_i\ge t\}=1/t^\lft$ for
$\ddl=\pm1$,
$i=1,\ldots,d$ and $t>0$. Let $\tilde\rdd$ be an excess measure on
$\rbb^d\sm\{\zerob\}$. Assume that for each open orthant $Q_\Delta$,
$\Delta\in\{-1,1\}^d$, the restrictions of $\tilde\rdd$ and $\rdd$ to
$Q_\Delta$ have the same univariate margins. One may then choose
$\tilde F$ such that its margins are tail asymptotic to $F_0$, such
that the sample clouds from $\tilde F$ converge to the Poisson point
process $\tilde N$ with mean measure $\tilde\rdd$, and such that the
sample clouds from the d.f. $\tilde G=\tilde F\circ K$ converge onto the
limit set $E_{\lft,\q}$ in \eqref{qE}.
\end{proposition}

\begin{pf}
We sketch the
construction. Choose $\hat F\in\DC_\lft(\tilde\rdd)$ with density
$\hat f$
such that the sample clouds from $\hat F$ scaled by $c_n$ converge to
$\tilde N$. For $\Delta\in\{-1,1\}^d$, let $U_\Delta$ be the image of
the union $U$ in $Q_\Delta$ by reflecting coordinates for which
$\Delta_j=-1$ (see Figure~\ref{fcase1} for an illustration). Let
$\tilde f$ agree with $\hat f$ on the $2^d$ sets $U_\Delta$ and with
$f$ elsewhere, so that, by the remark above on the convergence of the
scaled sets $U$, $\tilde f$ and $\hat f$ differ only on an
asymptotically negligible set. Alter $\tilde f$ on a bounded set to
make it a probability density. Then the sample clouds from $\tilde F$
scaled by $c_n$ converge to $\tilde N$. In the corresponding
partition $(A_n)$ on ${\mathbf x}$-space we only change the measure on the
``tiny'' blocks $(s_{n-1},s_{n+1})^d$ (with $s_{n-1}\sim s_{n+1}$)
around the positive diagonal, and their reflections. Hence, the scaled
sample clouds from $\tilde G=\tilde F\circ K$ converge onto
$E_{\lft,\q}$.
\end{pf}

%
%
\begin{figure}
\begin{tabular}{@{}cc@{}}

\includegraphics{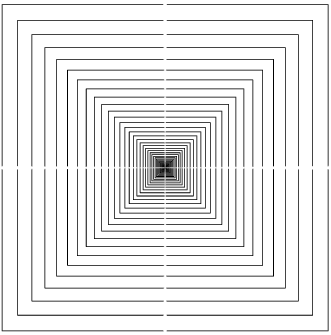}
&\includegraphics{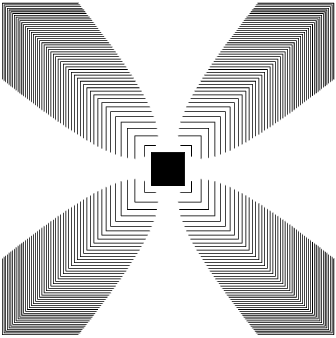}\\
(a)&(b)\\
\end{tabular}
\caption{Panel~(a): blocks
$[-t_n,t_n]^2=[-\mathrm{e}^{\sqrt{n}},\mathrm{e}^{\sqrt{n}}]^2$,
$n\in\{2,4,6,\ldots,100\}$ in ${\mathbf z}$-space with the sets
$[-\mathrm{e}^{\sqrt{n}}/n,\mathrm{e}^{\sqrt{n}}/n]\times[\mathrm{e}^{\sqrt{n-1}},\mathrm{e}^{\sqrt{n}}]$
and their reflections, which correspond to the blocks of partition
$(B_n)$ intersecting the axes, deleted. Panel~(b): image of the blocks
($[-\sqrt{n},\sqrt{n}]^2$) and deleted blocks ($[-\sqrt{n}+\log
n,\sqrt{n}-\log n]\times[\sqrt{n-1},\sqrt{n}]$) from Panel~(b) in
${\mathbf x}$-space under the meta transformation $K_0(s)=\mathrm{e}^s$, $s>s_0$.  (a) \textbf{z}-space (regular). (b) \textbf{x}-space.}
\label{fcase2} 
\end{figure}


We now discuss the second case: the block partition $(B_n)$ in
${\mathbf z}$-space is regular, but $(A_n)$ is not. Figure~\ref
{fcase2} depicts
sequences of cubes $s_nC$ and $t_nC$ in $\rbb^2$ on which partitions
$(A_n)$ and $(B_n)$ are based in the special case when $s_n=\sqrt{n}$
and $t_n=K_0(s_n)=\mathrm{e}^{s_n}=\mathrm{e}^{\sqrt{n}}$, along with subintervals
$[-\mathrm{e}^{\sqrt{n}}/n,\mathrm{e}^{\sqrt{n}}/n]\times\{\mathrm{e}^n\}$ in ${\mathbf
z}$-space mapping
onto $[-\sqrt{n}+\log n,\sqrt{n}-\log n]\times\{\sqrt{n}\}$ in
${\mathbf x}$-space, which correspond to the partition blocks
intersecting the
coordinate axes. We now have the following result.

In the standard set-up with d.f.s $F$ and $G=F\circ K$, margins $F_0$ and
$G_0=F_0\circ K_0$ and densities $f$ and $g$ there are two formulations
of the asymptotic behaviour, analytic and probabilistic. The analytic
formulation states that $f(r\wb)/f(r\eb)$ converges for $r\to\nf$
to a
limit function $h$ of the form $c/n_D^{\lft+d}$ where $n_D$ is a
continuous gauge function (of a bounded open star-shaped set $D$
containing the origin). The behaviour of the meta density is quite
different: $g(r\ub)/g(r\eb)$ converges to zero on the complement of the
compact set $E=E_{\lft,\q}$ and tends to $+\nf$ on the interior of this
set. The probabilistic version treats the asymptotic behaviour of
sample clouds. The sample clouds from $F$, scaled by $a_n$, converge to
a Poisson point process with intensity $c_0h$; the sample clouds from
$G$ scaled by $b_n$ converge onto $E$. Here
\begin{equation}\label{q0cn}
1-F_0(c_n)=1-G_0(b_n)=1/n,\qquad c_n=K_0(b_n).
\end{equation}
 The figures show
that the meta transformation $K$ moves points outside the coordinate
planes towards the diagonals. The block partitions introduced above
enable us to make this more precise.

\begin{proposition}\label{thc2} Let $F$ and $G=F\circ K$ satisfy the
standard set-up
with densities $f$ and~$g$. There exists a perturbation $\tilde F$ with
density $\tilde f$ and meta d.f. $\tilde G=\tilde F\circ K$ with density
$\tilde g$ with the properties:
\begin{itemize}
\item$\tilde f\le f$ outside the
unit ball, and hence $\tilde g\le g$ outside a bounded set;
\item
$\tilde f\equiv f$ on the $d$ coordinate planes, and hence $\tilde
g\equiv g$ on the coordinate planes;
\item for any unit vector $\z$
with nonzero components and any sequences ${\mathbf z}_n=r_n\z_n$ with
$r_n\to\nf$ and $\z_n\to\z$ eventually $\tilde f({\mathbf
z}_n)\equiv
f({\mathbf z}_n)$;
\item for any unit vector $\x$ with non-zero components
which does not lie on one of the $2^d$ diagonal rays and any sequence
${\mathbf x}_n=r_n\x_n$ with $r_n\to\nf$ and $\x_n\to\x$
eventually $\tilde
g({\mathbf x}_n)=0$;
\item the sample clouds from $\tilde F$ scaled by $c_n$
in~\textup{(\ref{q0cn})} converge to a Poisson point process with intensity
$c_0h$;
\item the sample clouds from $\tilde G$ scaled by $b_n$
in~\textup{(\ref{q0cn})} converge onto the diagonal cross $E_\times$
in~\textup{(\ref{qrsdx})}.
\end{itemize}
\end{proposition}

\begin{pf}
The construction is simple. We
introduce a regular partition for $F$, and delete the mass in the atoms
$B_n$ which intersect a coordinate plane $\{z_i=0\}$, replacing it on a
thin ridge around the coordinate planes to ensure that the second
condition holds. In the end, we increase the density on the unit ball
outside the coordinate planes to ensure that the new function $\tilde
f$ is a probability density. We shall now give the details for
dimension $d=2$. We focus on the positive horizontal axis. Choose
$t_n=n$. Choose $t_{n1}=r_n=\mathrm{o}(n)$ but so that $s_{n1}\sim
s_n=K_0\inv(t_n)$. This is possible since slow variation of
$L_0=K_0\inv$ implies that $L_0(\e t_n)/L_0(t_n)\to1$ for any $\e>0$.
Set $\tilde f\equiv0$ on the atom $(n-1,n)\times(-r_n,r_n)$ but keep
$\tilde f\equiv f$ on the rectangle $(n-1,n)\times(-\mathrm{e}^{-n},\mathrm{e}^{-n})$. Do
this for $n\ge n_0$, and do the same for the three other halfaxes. The
first four conditions hold by construction. The fifth condition follows
by the $\LB^1$ convergence of $\tilde h_n(\wb)=\tilde
f(r_n\wb)/f(r_n\eb)\to h(\wb)$ for $r_n\to\nf$ on $\e B^c$ by
Lebesgue's theorem on dominated convergence with variable bound
$h_n=f(r_n\wb)/f(r_n\eb)$. The last condition holds if we delete all
mass in the atoms $B_n$ intersecting the axes since that implies that
$\tilde g$ vanishes identically on the atoms
$[s_{n-1},s_n]\times[-s_{n1},s_{n1}]$ whose union for $n\ge n_0$ for
any $r>1$ and $\e>0$ eventually contains the sector
$S_{r\e}\ssdd\{x_1>0\}$ bounded by the halflines $x_2=\pm
\mathrm{e}^{-\e}(x_1-r)$, $x_1\ge r$. The mass in the rim between two successive
cubes, with the exception of the atoms intersecting the axes, is moved
towards the diagonals. This establishes the last statement. We leave it
to the reader to check that the subrectangles on which the density is
retained are so thin that they do not influence the asymptotic
behaviour of the sample clouds from the meta distribution. This is a
univariate issue. If $1-F_0^o(t)=\mathrm{O}(\mathrm{e}^{-t})$ for $t\to\nf$ then the meta
d.f. $G_0^o=F_0^o\circ K_0$ satisfies $n(1-G_0^o(\e b_n))\to0$ for any
$\e>0$. \end{pf}

The incompatibility of the partitions $(A_n)$ and $(B_n)=(K(A_n))$
introduced in this section gives one technical explanation for the
peculiar sensitivity of the limit shape of sample clouds from the meta
distribution. If we regard the atoms of the partition $(B_n)$ as nerve
cells, then regularity of $(A_n)$ will make the region around the
coordinate planes in ${\mathbf z}$-space far more sensitive than the remainder
of the space, and it is not surprising that cutting away these regions
has drastic effects on the limit.

\subsection{Mixtures}\label{sdom3}

In the standard set-up, there is a heavy-tailed d.f. $F$ with density $f$
and a light-tailed meta d.f. $G=F\circ K$ with density $g$. The margins
of $f$ are equal and symmetric, and so are the margins of~$g$. The
sample clouds from $F$ converge in distribution to a Poisson point
process on $\rbb^d\sm\{\zerob\}$ with a homogeneous mean measure; the
sample clouds from $G$ converge onto the compact set
$E_{\lft,\q}\ssdd[-1,1]^d$. By a perturbation of the density $f$, we may
obtain a new meta density whose sample clouds converge onto the
diagonal cross $E_\times$, see Theorem~\ref{thc2}, even though the
perturbation is so small that the asymptotics of the heavy-tailed
distribution are preserved. In this section, using a different
technique, we perform an additional perturbation. The asymptotics of
the heavy-tailed d.f. are again preserved, but the sample clouds from the
meta distribution now converge onto the union of the diagonal cross and
a compact star-shaped set $A$, the closure of an open star-shaped set
in $(-1,1)^d$ with a continuous gauge function.

Start with a pair of d.f.s $F$ and $G=F\circ K$ in the standard set-up.
Suppose $\hat F\in\DC_\lft(\rdd)$ has margins $\hat F_i$, which are tail
asymptotic to the margins $F_0$ of $F$ at $\pm\nf$ and which are
continuous and strictly increasing. The margins of the meta d.f. $\hat
G=\hat F\circ K$ have the same properties: continuous, strictly
increasing and tail asymptotic to $G_0$. We assume that the limit set
for $\hat G$ is the diagonal cross. We shall replace $\hat F$ by a d.f.
$\tilde F$, where $\mathrm{d}\tilde F=\mathrm{d}\hat F+\mathrm{d}F^o$ outside a bounded set and
where the marginal tails of $F^o$ are negligible compared to those of
$\hat F$, and hence so too the marginal tails of $G^o=F^o\circ K$ with
respect to those of $\hat G=\hat F\circ K$. It follows that $\tilde F$
and $\hat F$ have the same asymptotics, but it is possible to choose
$G^o$ to have the prescribed limit set $A$. The limit set of $\tilde G$
then is the union, $E_\times\cup A$.

\begin{example}\label{em1}
Let $G^o$ have density $g^o({\mathbf
x})=g_*(\|
{\mathbf x}\|_\nf)$
where $g_*$ is continuous, strictly decreasing and positive, and varies
rapidly. The margins of $G^o$ are equal, and have a symmetric density
$g_0^o(s)\sim(2s)^{d-1}g_*(s)$ because of the rapid variation of $g_*$.
Suppose $g_0(s)\sim as^b\mathrm{e}^{-ps^q}$ with $a,p,q$ positive. Let
$g_*(s)\sim g_0(s)/s^d$. Then $g_0^o(s)\ll g_0(s)$, but for $c>1$ we find
$g_0(cs)\ll g_0^o(s)$ since $\log s\ll s^q$ implies $\mathrm{e}^{-p(cs)^q}\le
s^2\mathrm{e}^{-ps^q}$ eventually. Let $g_0(c_n)\sim1/n$ and
$g^o_0(a_n)\sim1/n$. Then $a_n\le c_n\le cc_n$ eventually for any
$c>1$, and hence $a_n\sim c_n$. Let $1-G_0(b_n)=1/n$. Then $b_n\sim
c_n$ by univariate EVT. The sample clouds from the d.f. $\hat G$ above,
scaled by $a_n$ (or $b_n$ or $c_n$), converge onto the diagonal cross
$E_\times$, but the sample clouds from $G^o$ with the same scaling
converge onto the cube $[-1,1]^d$, even though the tails of the margins
of $G^o$ are negligible compared to those of $\hat G$. It follows that
the heavy-tailed d.f.s $\hat F$ and $\tilde F$ have the same asymptotics,
the margins are tail asymptotic and $\tilde F\in\DC_\lft(\rdd)$
with the
same scaling as $\hat F$. For the meta distributions $\hat G=\hat
F\circ K$ and $\tilde G=\tilde G\circ K$ with $\mathrm{d}\tilde G({\mathbf
x})=\mathrm{d}\hat
G({\mathbf x})+g_*(\|{\mathbf x}\|_\nf)\,\mathrm{d}{\mathbf x}$, the situation is
different. The sample
clouds from $\tilde G$ converge onto the union of the diagonal cross
and the coordinate cube.
\end{example}

If we replace the d.f. with density $g_*(\|{\mathbf x}\|_\nf)$ by a d.f.
which has
density $g_*(n_D({\mathbf x}))$ outside a bounded set, the sample clouds
converge onto the union of $E_\times$ and the closure of $D$.

\begin{proposition}\label{thmix} Let $F$ and $G=F\circ K$ with margins
$F_0$ and $G_0$
satisfy the standard set-up. Let $1-G_0(b_n)=1/n$. Let $\hat
F\in\DC_\lft(\rdd)$ have margins $\hat F_i$ which are continuous and
strictly increasing and tail asymptotic to $F_0$ at $\pm\nf$. Let
$\hat
G$ have limit set $E_\times$. Let $D\ssdd(-1,1)^d$ be an open star-shaped
set which contains the origin and has a continuous gauge function.
There exists a continuous strictly decreasing positive function $g_*$
on $[0,\nf)$ which varies rapidly such that
\begin{equation}\label{q3gg}
s^{d-1}g_*(s)/g_0(s)\to0,\qquad s^{d-1}g_*(s)/g_0(cs)\to\nf,\qquad c>1,
\end{equation}
where $g_0$ is the marginal density of $G$. Let $\tilde G$ be
a d.f.
such that $\mathrm{d}\tilde G({\mathbf x})=\mathrm{d}\hat G({\mathbf x})+g_*(n_D({\mathbf
x}))\,\mathrm{d}{\mathbf x}$ outside a
bounded set, and define $\tilde F=\tilde G\circ K\inv$. Then $\tilde F$
has the same asymptotics as $\hat F$ and the sample clouds from $\tilde
G$ scaled by $b_n$ converge onto the closure of $E_\times\cup D$.
\end{proposition}

\begin{pf}The tails of the margins of $F^o$ are
negligible with respect to
the tails of the margins of $\hat F$. This ensures that $\tilde F$ and
$\hat F$ have the same asymptotics. Outside a bounded set the sample
from $\tilde G$ is the superposition of a sample from $\hat G$ and from
$G^o$. Hence, the scaled sample clouds from $\tilde G$ converge onto
the union of $E_\times$ and the closure of $D$. It remains to find a
function $g_*$ which satisfies~(\ref{q3gg}). Choose $g_*(s)\sim
g_0(s)/s^d$. Then the first limit relation holds. Write $g_0=\mathrm{e}^{-\jfd}$.
By assumption, $\jfd$ varies regularly with exponent $\q>0$. Hence,
$\jfd(s)\ge s^{\q/2}$ eventually and for any $c>1$ there exists a
constant $s_0$ such that
\[
\jfd(cs)-\jfd(s)\ge(c^\q-1)s^{\q/2}/2\ge2\log s,\qquad s\ge s_0.
\]
This implies $g_0(cs)\ll
g_0(s)/s\sim s^{d-1}g_*(s)$.
\end{pf}

Although the shape of the limit set is rather unstable under even
slight perturbations of the original distribution, one may note the
persistence of the diagonal cross as a subset of the limit set. For
scaling constants $b_n$ in~(\ref{q0cn}), the univariate projections of
the sample clouds $M_n=\{\XB_1/b_n,\ldots,\XB_n/b_n\}$ converge onto
$[-1,1]$. Hence, the limit set $E$ for the multivariate sample clouds
$M_n$, if it exists, has univariate projections $[-1,1]$. One may use
the invariance principle for limit distributions in multivariate EVT to
show why the set $E$ often contains the diagonal cross $E_\times$. We
shall use the ideas expressed in Figure~\ref{fcase1}.

\begin{proposition}\label{pdiag} Let $F$ and $G=F\circ K$ with margins
$F_0$ and $G_0$
satisfy the conditions of the standard set-up. Let $\tilde F$ have
margins $F_i$ which are continuous and strictly increasing and tail
asymptotic to $F_0$. Assume the sample clouds from $\tilde F$ scaled by
$c_n$ converge to a Poisson point process with mean measure $\tilde
\rdd$,
where $\tilde\rdd$ charges $(0,\nf)^d$. If the sample clouds from
$\tilde
G=\tilde F\circ K$ can be scaled to converge onto a limit set $\tilde
E$ with coordinate projections $[-1,1]$, then $E$ contains the point
$\eb=(1,\ldots,1)$.
\end{proposition}

\begin{pf}
We make use of the block partitions of
Section~3.2. Consider the situation sketched in Figure 2. It suffices
to look at the positive orthant. Consider cubes
$C^A_n:=(s_n-Ma(s_n),s_n+Ma(s_n))^d$ centered at diagonal points
$s_n\eb$ for some $M>1$, where $a(s)$ is the scale function of the
marginal d.f. $G_0$. Recall that $a'(s)\to0$ and hence $a(s)/s\to0$ as
$s\to\nf$ which gives $C_n^A/s_n\to\{\eb\}$ for $s_n\to\nf$. The
corresponding points in ${\mathbf z}$-space are centered at the
diagonal points
$t_n\eb$ with $t_n=K_0(s_n)$, and given by
$C_n^B=K(C_n^A)=(K_0(s_n-Ma(s_n)),K_0(s_n+Ma(s_n)))^d$. Set
$t_n^\pm=K_0(s_n\pm Ma(s_n))$. Then
\[
\bigl(1-F_0(t_n^\pm)\bigr)/\bigl(1-F_0(t_n)\bigr)=\bigl(1-G_0\bigl(s_n\pm Ma(s_n)\bigr)\bigr)/\bigl(1-G_0(s_n)\bigr)\to
\mathrm{e}^{\mp M}
\]
by univariate EVT applied to $G_0$. Regular variation of $1-F_0$ then
gives $t_n^\pm/t_n\to \mathrm{e}^{\pm M/\lft}$. Thus
$C^B_n/t_n\to(\mathrm{e}^{-M/\lft},\mathrm{e}^{M/\lft})^d$. For large $M$, this limit cube
constitutes a large part of $(0,\nf)^d$ in ${\mathbf z}$-space. Take $A>1$
large and $M/\lft>A$. Let $t_n=c_n$ and $s_n=b_n$ be the scaling
constants defined by $1-G_0(s_n)=1-F_0(t_n)=1/n$. Let $\e>0$ be small
and $n$ so large that $a(s_n)/s_n<\e$. Then the sample points
$\ZB_1/t_n,\ldots,\ZB_n/t_n$ in the cube $(\mathrm{e}^{-A},\mathrm{e}^A)^d$ yield points
$\XB_i/s_n$ in $\eb+(-\e,\e)^d$. Hence, $\eb\in E$ if $\tilde\rdd$
charges $(0,\nf)^d$.
\end{pf}

\section{Discussion}\label{sconc}

In situations where chance plays a role the asymptotic description
often consists of two parts, a deterministic term, catching the main
effect, and a stochastic term, describing the random fluctuations
around the deterministic part. Thus, the average of the first $n$
observations converges to the expectation; under additional assumptions
the difference between the average and the expectation, blown up by a
factor $\sqrt n$, is asymptotically normal. Empirical d.f.s converge to
the true d.f.; the fluctuations are modeled by a time-changed Brownian
bridge. For a positive random variable, the $n$-point sample clouds
$N_n$, scaled by the $1-1/n$ quantile, converge onto the interval
$[0,1]$ if the tail of the d.f. is rapidly varying; if the tail is
asymptotic to a von Mises function then there is a limiting Poisson
point process with intensity $\mathrm{e}^{-s}$.

Convergence to the first-order deterministic term in these situations
is a much more robust affair than convergence of the random
fluctuations around this term. So it is surprising that for meta
distributions perturbations of the original distribution which do not
affect the second-order fluctuations of the sample cloud at the
vertices may drastically alter the shape of the limit set, the
deterministic first-order term. This paper tries to cast some light on
the sensitivity of the meta distribution and the limit set $E$ to small
perturbations of the original distribution.

Bivariate asymptotics are well expressed in terms of polar coordinates.
Two points far off are close together if the angular parts are close
and if the quotient of the radial parts is close to one. This geometry
is respected by certain partitions. A partition is regular if points in
the same atom are uniformly close as one moves out to infinity. Call
probability distributions \textit{equivalent} if they give the same or
asymptotically the same weight to the atoms of a regular partition.
Equivalent distributions have the same asymptotic behaviour with
respect to scaling.

This paper compares the asymptotic behaviour of a heavy-tailed density
with the asymptotic behaviour of the associated\vadjust{\goodbreak} meta density with
light-tailed margins. Small changes in the heavy-tailed density,
changes which have no influence on its asymptotic behaviour, may lead
to significant changes in the asymptotic behaviour of the meta
distribution. We show that regular partitions for the heavy-tailed
distribution and for the light-tailed meta distribution are
incommensurate. The atoms at the diagonals in the light-tailed
distribution fill up the quadrants for the heavy-tailed distribution;
atoms at the axes in the heavy-tailed distribution fill up the four
segments between the diagonals for the light-tailed distributions.
Section~\ref{sdom2} shows how equivalent distributions in one space
give rise to different asymptotic behaviour in the other.

In our approach, the asymptotic behaviour in both spaces is
investigated by rescaling. In the heavy-tailed world one obtains a
limiting Poisson point process whose mean measure is an excess measure
$\rdd$ which is finite outside centered disks in the plane; in the
light-tailed world the sample clouds converge onto a star-shaped limit
set $E$. In the standard set-up the only relation between $\rdd$ and $E$
is the parameter $\lft$. This parameter describes the rate of decrease of
the heavy-tailed marginal distributions; it also is one of the two
parameters which determine the shape of the limit set $E$.

We can offer two explanations for the incompatibility of the
asymptotics of a heavy-tailed d.f. $F$ and the light-tailed meta d.f. $G$.
\begin{longlist}
\item[(1)] The geometric explanation is that the meta transformation does not
preserve direction. A~ray in ${\mathbf x}$-space which does not lie in a
diagonal plane is transformed into a curve whose direction is
asymptotic to a halfaxis. Conversely a ray in ${\mathbf z}$-space which does
not lie in a coordinate plane lies in one of the $2^d$ open orthants
$Q_\DG$ and is transformed by $K\inv$ into a curve which is asymptotic
to the diagonal ray in the center of the orthant. This geometric
distortion also occurs if one moves from heavy tails to less heavy
tails, increasing the parameter $\lft$ of regular variation, but to a
lesser extent. See~\cite{fest}.
\item[(2)] The probabilistic interpretation is that the limit set $E$ describes
the intermediate extremes whereas the limiting Poisson point process
$N$ describes the asymptotics of the extreme order statistics. The d.f.
$F^o$ in Section~\ref{sdom3} contributes to the intermediate order statistics,
but not to the extremes.
\end{longlist}

\section*{Acknowledgements} We would like to thank the three referees
and the editor for their valuable comments on the earlier version of
the paper. The authors acknowledge financial support of the Institute
for Mathematical Research (FIM) and RiskLab, Switzerland. Paul
Embrechts, as Senior SFI Professor, would like to thank the Swiss
Finance Institute for financial support.

%

\printhistory

\end{document}